\documentclass[sn-mathphys]{sn-jnl}

\usepackage{physics} 
\usepackage{mathtools} 
\usepackage[english]{cleveref} 
\usepackage{enumitem}

\jyear{2021}%

\theoremstyle{thmstyleone}%
\newtheorem{theorem}{Theorem}
\newtheorem{lem}[theorem]{Lemma}%

\theoremstyle{thmstyletwo}%
\newtheorem{remark}{Remark}%

\theoremstyle{thmstylethree}%
\newtheorem{definition}{Definition}%

\newcommand{\J}{\mathcal{J}} 
\newcommand{\D}{\mathrm{D}} 
\newcommand{\dif}{\,\mathrm{d}} 
\DeclareMathOperator{\sing}{Sing} 
\newcommand{\cals}{\mathcal{S}} 
\newcommand{\mfk}{\mathbb{R}^3\times\mathcal{S}} 
\newcommand{\moddeg}{\((\hspace{-7pt}\mod 2)\)-degree } 
\newcommand{\cmj}{C_{m}^{j}} 
\newcommand{\sball}{B_{\frac{\alpha}{2}}^2\times\{0\}} %
\newcommand{\ball}{B_{\alpha}^2\times\{0\}}
\newcommand{\ring}{\left(B_{\alpha}^2\setminus B_{\frac{\alpha}{2}}^2\right)\times\{0\}}

\begin{document}

\title{On prescribing the number of singular points in a Cosserat-elastic solid}

\author{\fnm{Vanessa} \sur{Hüsken}\footnote{E-mail: vanessa.huesken@uni-due.de}}

\affil{\orgdiv{Faculty of Mathematics}, \orgname{University Duisburg-Essen}, \orgaddress{\street{Thea-Leymann-Str. 9},  \postcode{45127}, \city{Essen}, \country{Germany}}}


\abstract{In a geometrically non-linear Cosserat model for micro-polar elastic solids, we insert dipole pairs of singularities into smooth maps and control the amount of Cosserat energy needed to do so. We use this method to force an arbitrary number of singular points into Cosserat-elastic solids by prescribing smooth boundary data. Throughout this paper, we exploit connections between harmonic maps and Cosserat-elastic solids, so that we are able to adapt and incorporate ideas of R. Hardt and F.-H. Lin, as well as of F. Béthuel.}



\maketitle


\section{Introduction and statement of results.}\label{sec1:introduction}

Cosserat elasticity is a well known class of models in elasticity theory, whose foundations were laid at the beginning of the \(20^{th}\) century by the Cosserat brothers. The geometrically non-linear model for micro-polar elastic solids being discussed in this paper is a type of Cosserat elasticity that has first been studied in the context of calculus of variations by P. Neff, for example in \cite{neff1}. Its basic concept is the following. An elastic body in its original state is described as a subset \(\Omega\subset\mathbb{R}^3\). It can be deformed by shifting each point \(x\in\Omega\) to its new location \(\varphi(x)\in\mathbb{R}^3\). Moreover, the micro-polar structure of the body allows each point to undergo some micro-rotation (without deforming the body any further), meaning that to each point \(x\), there is attached an orthonormal frame, which is free to rotate by an orthogonal matrix \(R(x)\in SO(3)\). The micro-rotation being in \(SO(3)\), rather than using infinitesimal rotations in the corresponding Lie-algebra of skew-symmetric matrices, ultimately leads to the geometric non-linearities in the Euler-Lagrange equations of the model. Both deformation and micro-rotation cause material stresses, measured in terms of \(R^T\cdot\D\varphi -I_3\) and \(R^T\cdot\D R\), respectively. Leaving additional external forces and moments aside (as it was discussed in \cite{gastel2019regularity}), summing up the energy stored in the body, the \emph{Cosserat energy functional} is given by \[\J_{\Omega}(\varphi,R) = \norm{P(R^T\cdot\D\varphi-I_3)}_{L^2(\Omega)}^2 + \lambda\norm{R^T\cdot\D R}_{L^p(\Omega)}^p,\] with constant \(\lambda>0\), parameter \(p\geq 2\) and linear operator \(P\colon\mathbb{R}^{3\times 3}\to\mathbb{R}^{3\times 3}\), describing a weighted sum of the deviatoric symmetric part and the skew-symmetric part of a matrix as well as a diagonal matrix of its trace: \[P(A)= \sqrt{\mu_1}\operatorname{dev sym}(A) + \sqrt{\mu_c}\operatorname{skew}(A) + \frac{\sqrt{\mu_2}}{3}\tr(A)\cdot I_3,\] with positive material constants \(\mu_1,\mu_c,\mu_2>0\).\\

The existence of minimizers of this Cosserat energy on a bounded Lipschitz domain \(\Omega\subset\mathbb{R}^3\) was proven in \cite{neff2006existence}. Further aspects of the model itself and the existence of Cosserat energy minimizers are also discussed in \cite{nbo}.

When studying regularity of minimizers, Gastel recently observed a connection between the Cosserat problem and \(p\)-harmonic maps, which is a well studied area in Geometric Analysis. In the case \(p=2\) (\(\lambda=1\) without loss of generality), when all constants are assumed equal (\(\mu_1=\mu_c=\mu_2\)), he found the following (cf. \cite{gastel2019regularity}): On one hand, he showed Hölder-continuity for all minimizers on the whole domain \(\Omega\). On the other hand, he gave an example of a critical point (meaning a weak solution of the Euler-Lagrange equations) of the Cosserat energy for \(\Omega=B^3\), \(p=2\) and \(\mu_1=\mu_c=\mu_2=1\), whose micro-rotational part exhibits a point singularity at the origin. So in contrast to minimizers, regularity of critical points should be an issue.\\

Note that with this particular choice of constants \(P(\cdot)\) becomes the identity and \[\J_{\Omega}(\varphi,R) = \int\limits_{\Omega} \abs{R^T\cdot\D\varphi - I_3}^2 + \abs{\D R}^2 \dif x.\] In Geometric Analysis, many results are known about (non-)regularity of harmonic mappings (i.e. weak solutions for the Euler-Lagrange equations of the Dirichlet integral). Having in mind several of them, concerning harmonic mappings into the standard sphere \(S^2\), the starting point of our research is the question: How "big" can the singular set \(\sing(f)\) of a critical point \(f=(\varphi,R)\) of the Cosserat energy get? (In this situation, \(\sing(f)\) denotes the set of points, where \(f\) fails to be locally in \(C^{1,\mu}\times C^{0,\mu}\) for any \(\mu\in(0,1)\), its elements are called singularities. Similarly, \(\sing(R)\) denotes the set, where \(R\) fails to be locally in \(C^{0,\mu}\) for any \(\mu\in(0,1)\).) \smallskip

An idea for being able to use the vast machinery of results about the regularity of harmonic mappings into \(S^2\) is to observe a connection between \(S^2\) and the set \[\cals \coloneqq \{A\in SO(3) : A \text{ describes a }180^{\circ}\text{-rotation around some axis in }\mathbb{R}^3\}.\] By identifying each rotation in \(\cals\subset SO(3)\) with its axis of rotation, we obtain a two-fold covering of the manifold \(\cals\), given by \(F\colon S^2 \to \cals\), \(q\mapsto 2q\otimes q - I_3\). A quick calculation in local coordinates shows that \(F\) is locally isometric up to the factor \(\sqrt{8}\). Moreover, a well known fact from Algebraic Topology implies that, if the domain \(\Omega\) is simply connected and locally path-connected, any continuous mapping can be lifted \cite[][Thm. 6.1 \& Cor. 6.4, p. 26f.]{greenberg2018algebraic}. To be precise, for the covering \(F\) and any continuous mapping \(R\colon\Omega\to\cals\), there exist two continuous mappings \(\eta_{1,2}\colon\Omega\to S^2\), \(\eta_1=-\eta_2\), such that \(R=F\circ\eta_{1,2}\), as long as \(\Omega\) is simply connected and locally path-connected.

So instead of looking at the full variational Cosserat problem 
	\begin{equation}\label{eq:fullproblem}
		\J_{\Omega}(\varphi,R)=\int\limits_{\Omega} \abs{R^T\cdot\D\varphi-I_3}^2 +\abs{\D R}^2 \dif x 	\longrightarrow\min\quad \text{in } H^1(\Omega,\mathbb{R}^3\times SO(3)), \tag{$\mathcal{P}$}
	\end{equation}
we mostly work with the \emph{restricted Cosserat problem}
	\begin{equation}\label{eq: restricted problem}
		\J_{\Omega}(\varphi,R)=\int\limits_{\Omega} \abs{R^T\cdot\D\varphi-I_3}^2 +\abs{\D R}^2 \dif x 	\longrightarrow\min\quad \text{in } H^1(\Omega,\mfk). \tag{$\mathcal{P}^{\ast}$}
	\end{equation} 
Often, restricting a variational problem to a submanifold changes the Euler-Lagrange equations and thus is not a suitable method for finding results for the general problem. But \(\cals\) is a totally geodesic submanifold of \(SO(3)\). This fact implies (just like it is proven for harmonic mappings), that \emph{restricted minimizers} (i.e. minimizers of the restricted Cosserat problem (\ref{eq: restricted problem})) are at least still critical points of the full Cosserat problem (\ref{eq:fullproblem}). In general, they are not minimizers of (\ref{eq:fullproblem}).

In \cite{gastel2019regularity}, Gastel showed that the (interior) singular set of a Cosserat energy minimizer of the full problem (\ref{eq:fullproblem}) is a discrete set and in fact empty. But the line of reasoning made there to show discreteness holds true for restricted minimizers. With similar arguments, following the suggestions from \cite{steffen1991harmonicintroduction}, based on \cite{schoen1983boundary} in the context of harmonic maps, one can also show discreetness of the singular set at the boundary. So we expect only point singularities for restricted minimizers and contrary to \cite{gastel2019regularity}, in analogy to a result from \cite{hardt1986remark} for harmonic maps \(u\colon B^3\to S^2\), we derive the following statement. It shows that critical points of the Cosserat energy can be forced to have an arbitrary large number of singularities, by prescribing suitable smooth boundary data.\smallskip

	\begin{theorem}\label{thm:HL}
		For every number \(N\in\mathbb{N}\) there exist smooth boundary data \(g_0=\left(\varphi_0,R_0\right)\in C^{\infty}(\partial B^3,\mfk)\) with \(\deg(R_0)=0\), such that each (restricted) minimizer \(f=\left(\varphi,R\right)\) of the Cosserat energy \(\J_{B^3}(\cdot)\) in the class \(H^1_{g_0}(B^3,\mfk)\coloneqq\{g\in H^1(B^3,\mfk):\,g_{\vert\partial B^3}=g_0\}\)
		must have at least \(N\) singularities in its micro-rotational part \(R\).
	\end{theorem}

	\begin{remark} 
		The property \(\deg(R_0)=0\) emphasizes that the singularities, which we are about to enforce, do not appear simply due to elementary topological reasons, see the discussion in \cite[p. 15]{brezis1989s} for example, in regard to harmonic mappings \(u\colon\Omega\to S^2\). But, as \(\cals \simeq\mathbb{R} P^2\) is a non-orientable manifold, the concept of the classical Brouwer-degree \(\deg(\psi)\) of a mapping \(\psi\) between orientable manifolds, which is used for the deformation component \(\varphi\), needs to be modified for the micro-rotational component \(R\). Inspired by observations in \cite{olum1953mappings} and \cite[][§ 4]{milnor1997topology}, we define the \moddeg as follows.
	\end{remark}

	\begin{definition}\label{def:moddeg}
		Let \(\Omega\subset \mathbb{R}^3\) be a bounded, simply connected and locally path-connected set and let \(f=(\varphi,R)\colon\Omega\to\mfk\).
			\begin{enumerate}[label=(\roman*)]
				\item For \(R\in C^0(\Omega,\cals)\), there exists a lift \(n\colon\Omega\to S^2\), which means \(F\circ n=R\). Then the \emph{\moddeg} of \(R\) is given by \[\deg(R)\coloneqq \deg(n)\mod 2.\]
					\item For an isolated singularity \(a\in\sing(R)\), we define \[\deg_a(R)\coloneqq\deg(R_{\vert S^2_r(a)}) = \deg_a(n_{(a)}) \mod 2,\] where \(S^2_r(a)\subset\Omega\) is an arbitrary sphere of radius \(r>0\) around \(a\), such that the corresponding ball \(\overline{B^3_r}(a)\) does not contain any other singularities of \(f\), and \(n_{(a)}\) denotes the lift of \(R\) existing on \(\overline{B^3_r}(a)\setminus\{a\}\). 
			\end{enumerate}
	\end{definition}
	In both cases, the \moddeg lies in \(\mathbb{Z}/2\mathbb{Z}\). This definition has the advantage, that nice properties of the classical Brouwer-degree, like additivity and homotopy invariance, continue to hold.\\

Because we are going to use the concept of dipoles a lot throughout this paper, we also have to modify the original definition of a dipole as introduced in \cite{bcl1986harmonic} to fit into the situation of (restricted) Cosserat solids.
	\begin{definition}\label{def:dipol}
		Let \(\Omega\subset\mathbb{R}^3\) and \(f=(\varphi,R)\colon\Omega\to\mfk\) be as in \Cref{def:moddeg}. A pair of singularities \((P,N)\) of \(R\) is called a \emph{dipole} for \(R\), if there is an open bounded cylinder \(Z^3_r(q)\subset\Omega\), rotationally symmetric (of radius \(r>0\)) around the line segment \([P,N]\), such that
			\begin{enumerate}[label=(\roman*)]
				\item \([P,N]\subset Z^3_r(q)\) and \(Z^3_r(q)\) is centred at the centre \(q\) of \([P,N]\) ,
				\item \(Z^3_r(q)\) does not contain any further singularities of \(f\),
				\item \(\deg_P(R)=1=\deg_N(R)\) and the lift \(n_{(q)}\) of \(R\) (existing on \(\overline{Z^3_r}(q)\setminus\{P,N\}\)) has a classical dipole \((P,N)\), i.e. \(\deg_P(n_{(q)})=d=-\deg_N(n_{(q)})\) for a \(d\in\mathbb{Z}\).
			\end{enumerate}
	\end{definition}
A central method to prove \Cref{thm:HL} in \cref{sec:resmin} is inserting dipoles into a given smooth mapping while controlling the energy needed to do so. The details are stated in the following theorem, which will be proven in \cref{sec:dipoles}.

	\begin{theorem}\label{thm:dipoles}
		Let \(\Omega\subset\mathbb{R}^3\) be a bounded, simply connected and locally path-connected set. Let \(P,N\) be two distinct points in \(\Omega\), such that the line segment  \([P,N]\) lies fully in \(\Omega\). Then, for any mapping \(f=(\varphi,R) \in C^{\infty}(\Omega,\mathbb{R}^3\times\cals)\), there exists a sequence of mappings \[f_m=(\varphi_m,R_m)\in H^1(\Omega,\mathbb{R}^3\times\cals)\cap C^{\infty}(\Omega\setminus\{P,N\},\mathbb{R}^3\times\cals)\] with \((P,N)\) being a dipole for each \(R_m\), i.e. in particular it holds
			\begin{equation*}
				\deg_P(R_m) = 1 = \deg_N(R_m).
			\end{equation*}
		All mappings \(f_m\) agree with \(f\) outside of a neighbourhood \(K_m\) of \([P,N]\) that fulfils
			\begin{equation*}
				K_m \longrightarrow [P,N], \,m\to \infty,\quad \text{in Hausdorff-distance.}
			\end{equation*}
		Moreover,
			\begin{equation*}
				\lim\limits_{m\to\infty} \J_{\Omega}(f_m) \leq \J_{\Omega}(f)+64\pi\abs{P-N}.
			\end{equation*}
	\end{theorem}


\section{Construction of dipoles}\label{sec:dipoles}
As mentioned above, a key ingredient in the construction of suitable boundary data for the proof of \Cref{thm:HL} is the insertion of dipole pairs of singularities, each with \moddeg 1, into smooth maps. \Cref{thm:dipoles} gives us a tool for doing so while using a controlled amount of Cosserat energy, depending only on the dipole's length. The main part of this paper will consist of its proof, as it contains some technical intricacies.

\begin{proof}[Proof of \Cref{thm:dipoles}.]
	This proof is divided into three steps: First, we present a construction, that was used by F. Béthuel in \cite{BETHUEL1990269} to remove a dipole from a given map with a controlled amount of (Dirichlet) energy. Working in the other direction, it gives rise to a sequence of Lipschitz mappings with the desired singularities of \moddeg 1 inserted. Additionally, the mappings of the sequence exhibit further singularities of degree 0. Second, we calculate the estimates for the Cosserat energy needed. During the last step, we use some approximation results from \cite{bethuelapprox} to replace each Lipschitz mapping of the sequence with an approximation in order to get rid of the additional singularities of degree 0 and to gain the desired smoothness (except in \(P,N\)) without affecting the Cosserat energy.\\
	
	\textbf{Step 1}\label{par:stepone} (Construction). In \cite[][Lemma 2]{BETHUEL1990269}, F. Béthuel uses a cuboid construction together with a cube lemma \cite[][Lemma 3]{BETHUEL1990269} plus some calculations from the two-dimensional case in \cite{brezis1983large}. We can use exactly the same cuboid construction, together with the following modified cube lemma for the Cosserat energy which itself will be proved after having completed the proof of \Cref{thm:dipoles}.
	
		\begin{lem}[Cube Lemma]\label{lem:cube}		
			For \(\nu>0\), let \(C_{\nu} = [-\nu,\nu]^2\times[-2\nu,0]\). Consider a Lipschitz mapping \(f=(\varphi,R)\colon\partial C_{\nu}\to\mathbb{R}^3\times\mathcal{S}\) with \(\deg(R)=d_0\), \(d_0\in\mathbb{Z}/2\mathbb{Z}\).
		
			Then, for each \(\varepsilon>0\), there exists a constant \(\alpha_0\in(0,\nu)\), such that for any \(0<\alpha<\alpha_0\), there exists a Lipschitz mapping \(f_{\alpha}=(\varphi,R_{\alpha})\colon\partial C_{\nu}\to\mfk\) with
				\begin{align}
					\deg(R_{\alpha}) &= d_0+1 \mod 2,\nonumber\\
						f_{\alpha}=(\varphi,R_{\alpha}) &= (\varphi,R)=f \text{ in } \partial C_{\nu}\setminus\left(\ball\right)\nonumber\\
					\shortintertext{and}
						\int\limits_{\ball} 2\abs{R_{\alpha}^{T}\cdot \D\varphi - I_3}^2 &+ \abs{\D R_{\alpha}}^2 \dif\mathcal{H}^2 < 64\pi + \varepsilon.\label{eq:O-1-3-3}
				\end{align}
			Moreover, on \(\ring\) we have
				\begin{equation}\label{eq:O-1-3-4}
					\abs{\D R_{\alpha}}\leq const,
				\end{equation}
			and on \(\sball\)
				\begin{equation}\label{eq:O-1-3-5}
					\abs{\D R_{\alpha}(x,y,0)}^2 = \frac{64\alpha^4}{(\alpha^4+x^2+y^2)^2}.
				\end{equation}
		\end{lem}
		
	\noindent Following the notation from \cite{BETHUEL1990269} for the cuboid construction, with 
		\begin{itemize}
			\item \(d\coloneqq \abs{P-N}\),
				\item \(a_m\coloneqq\frac{d}{2(m-1)}\), \(m\in\mathbb{N}_{>1}\),
			\item \(K_m\) the cuboid around \([P,N]\): \(K_m\coloneqq \left[-a_m,a_m\right]^2\times\left[-a_m,d+a_m\right] \),
				\item \(K_m\) divided into \(m\) cubes: \(\cmj\coloneqq\left[-a_m,a_m\right]^2\times\left[(-1+2j)a_m,(1+2j)a_m\right],\)\\ \(j=0,\dots,m-1\),
			\item \(c_j\) the barycentre of \(\cmj\) and
				\item \(\pi_{m}^j\colon \cmj \to \partial\cmj\) the radial retraction with centre \(c_j\), given by\\
				\(\pi_{m}^j(x)= \frac{x-c_j}{\abs{x-c_j}_{\infty}}\cdot a_m + c_j\), \(\qquad\abs{x-c_j}_{\infty}=\max_{i=1,2,3}(\abs{x_i-(c_j)_i})\),
		\end{itemize}
	we iteratively use \Cref{lem:cube} on the single cubes \(\cmj\) implying that for each \(m\geq m_0\gg1\), there exists a sequence of Lipschitz mappings \((\tilde{f}_{m,\alpha})_{\alpha}\in H^1(\Omega,\mfk)\) (with \(\alpha\searrow 0\)) given by \[\tilde{f}_{m,\alpha} = (\tilde{\varphi}_{m,\alpha},\tilde{R}_{m,\alpha}) \coloneqq \begin{cases}
		(\varphi,R) = f, & \text{ in }\Omega\setminus K_m,\\
		f_{m,\alpha}\circ\pi_{m}^{j} , & \text{ in } K_m,
	\end{cases} \]
	where \[f_{m,\alpha}\colon\bigcup\limits_{j=0}^{m-1} \partial \cmj \to\mfk, \qquad f_{m,\alpha}=(\varphi,R_{m,\alpha}), \] with \(\sing(\tilde{f}_{m,\alpha})=\{P=c_0, c_1,\dotsc,c_{m-2},c_{m-1}=N\}\) and \vspace{-5pt}
		\begin{align*}
			\deg_P(\tilde{R}_{m,\alpha}) &=  1 = \deg_N(\tilde{R}_{m,\alpha}),\\
				\deg_{c_j}(\tilde{R}_{m,\alpha}) &= 0  \qquad\text{ for } j=1,\dotsc,m-2,\\
			\deg_{c_j}(\tilde{\varphi}_{m,\alpha}) &= 0  \qquad\text{ for } j=0,\dotsc,m-1,\\
				\shortintertext{as well as}
			\tilde{f}_{m,\alpha} &=f \qquad \text{ on } \partial K_m,
		\end{align*}
	due to the fact that during the construction, changes of the original mapping only happen on little discs (of radius \(\alpha\)) on the upper faces of the lower \(m-1\) cubes \(C_m^0,\dotsc,C_m^{m-2}\), so that the values on \(\partial K_m\) remain unaffected.
		
	\textbf{Step 2} (Calculation of Cosserat energy cost). Many of the calculations in this step follow ideas and estimates carried out in \cite{tarp} in the context of removing dipoles from given maps \(u\colon\Omega\to S^2\). Similarly, in order to calculate the Cosserat energy of \(\tilde{f}_{m,\alpha}\) on \(K_m\), we divide each cube \(\cmj\) into disjoint sets
		\begin{align*}
			B^3_{a_m}(c_j), & \\
				A_{m}^{j} &= \left( \cmj\setminus B^3_{a_m}(c_j) \right) \cap (\pi_{m}^{j})^{-1} \left( B^2_{\frac{\alpha}{2}}\times\{(-1+2j)a_m\} \right),\\
			D_{m}^{j} &= \left( \cmj\setminus B^3_{a_m}(c_j) \right) \cap (\pi_{m}^{j})^{-1} \left( B^2_{\frac{\alpha}{2}}\times\{(1+2j)a_m\} \right),\\
				E_{m}^{j} &= \left( \cmj\setminus B^3_{a_m}(c_j) \right) \cap (\pi_{m}^{j})^{-1} \left( \left(B^2_{\alpha}\setminus B^2_{\frac{\alpha}{2}}\right)\times\{(-1+2j)a_m\} \right),\\
			F_{m}^{j} &= \left( \cmj\setminus B^3_{a_m}(c_j) \right) \cap (\pi_{m}^{j})^{-1} \left( \left(B^2_{\alpha}\setminus B^2_{\frac{\alpha}{2}}\right)\times\{(1+2j)a_m\} \right)\\
				\shortintertext{and the rest}
			G_{m}^{j} &= \cmj\setminus\left(B^3_{a_m}(c_j)\cup A_{m}^{j}\cup D_{m}^{j}\cup E_{m}^{j}\cup F_{m}^{j}\right).
		\end{align*} 
	
	\noindent Note that different constants appearing in the following estimates are always denoted by the same \(\gamma\in\mathbb{R}_{\geq 0}\). For the Cosserat energy in \(B^3_{a_m}(c_j)\), we get
		\begin{align}
			&\J_{B^3_{a_m}(c_j)}(\tilde{f}_{m,\alpha}) = \int\limits_{B^3_{a_m}(c_j)} \hspace{-10pt} \abs{(R_{m,\alpha}\circ\pi_{m}^{j})^T\cdot \D(\varphi\circ\pi_{m}^{j}) - I_3}^2 + \abs{\D(R_{m,\alpha}\circ\pi_{m}^{j})}^2 \dif x\nonumber\\[-10pt]
				&= \int\limits_0^{a_m}\int\limits_{S^2_{\rho}(c_j)} \abs{(R_{m,\alpha}\circ\pi_{m}^{j})^T\cdot \D(\varphi\circ\pi_{m}^{j}) - I_3}^2 + \abs{\D(R_{m,\alpha}\circ\pi_{m}^{j})}^2\dif\mathcal{H}^2\dif\rho\nonumber\\[-10pt]
			&\leq \gamma a_m^3 + \operatorname{Lip}^2\left(\pi^{j}_{m \vert S^2_{\rho}(c_j)}\right) \int\limits_0^{a_m}\int\limits_{\partial\cmj} \Big\{2\abs{R_{m,\alpha}^T \cdot \D\varphi - I_3}^2\nonumber\\[-10pt]
				&\qquad\qquad\qquad\qquad\qquad\qquad\qquad\qquad + \abs{\D R_{m,\alpha}}^2\Big\} \operatorname{Jac}^{-1} \left(\pi^{j}_{m \vert S^2_{\rho}(c_j)}\right) \dif \mathcal{H}^2 \dif \rho\nonumber\\[10pt]
			&\leq \gamma a_m^3 + 9a_m\cdot 2\cdot \int\limits_{\partial\cmj\setminus\left(B^2_{\alpha}\times\{(-1+2j)a_m,(1+2j)a_m\}\right)} \hspace{-55pt}\abs{R^T\cdot \D\varphi-I_3}^2 + \abs{\D R}^2 \dif\mathcal{H}^2\nonumber\\
				&\qquad+ a_m\hspace{-5pt}\int\limits_{B^2_{\alpha}\times\{(-1+2j)a_m,(1+2j)a_m\}} \hspace{-45pt}\Big(1+\tfrac{\abs{y}^2}{a_m^2}\Big)^2\Big\{2\abs{R_{m,\alpha}^T\cdot\D\varphi - I_3}^2 + \abs{\D R_{m,\alpha}}^2\Big\} \dif \mathcal{H}^2(y)\nonumber\\[10pt]
			&\leq \gamma a_m^3 +18a_m\int\limits_{\partial\cmj} \abs{R^T\cdot \D\varphi-I_3}^2 + \abs{\D R}^2\dif\mathcal{H}^2 + 2a_m\left(1+\frac{\alpha^2}{a_m^2}\right)^2(64\pi+\varepsilon)\nonumber\\
				&\leq \gamma a_m^3 + 2a_m\left(1+\frac{\alpha^2}{a_m^2}\right)^2(64\pi+\varepsilon), \label{eq:O-1-4-13}
		\end{align}
	because of (\ref{eq:O-1-3-3}) and because the original \(f=(\varphi,R)\) is smooth in all of \(\Omega\). Therefore, \(\abs{R^T\cdot \D\varphi - I}^2\) and \(\abs{\D R}^2\) are bounded on \(K_{m}\), thus for each  \(j=0,\dotsc,m-1\)		
		\begin{equation*}
			\int\limits_{\partial\cmj} \abs{R^T\cdot 	\D\varphi-I_3}^2 + \abs{\D R}^2 \dif\mathcal{H}^2\leq \gamma\cdot a_m^2.
		\end{equation*}
		
	In \( G_{m}^{j}\) we have \(\tilde{f}_{m,\alpha} (x) =(f\circ\pi_{m}^{j})(x)\). Hence,
		\begin{equation}
			\J_{G_{m}^{j}}(\tilde{f}_{m,\alpha}) = \int\limits_{G_{m}^{j}} \abs{(R\circ\pi_{m}^j)^T\cdot \D(\varphi\circ\pi_{m}^j) - I_3}^2 + \abs{\D(R\circ\pi_{m,j})}^2 \dif x\leq \gamma\cdot a_m^3,\label{eq:O-1-4-14}
		\end{equation}
	as \((R\circ\pi_{m}^j)\in\mathcal{S}\subset SO(3)\), \(\D\varphi\) and \(I_3\), as well as \(\D R\) are bounded and \(\pi_{m}^j\) is Lipschitz in \(\cmj\setminus B^3_{a_m}\). Similarly we have
		\begin{alignat*}{2}
			\tilde{f}_{m,\alpha} &= f\circ\pi_{m}^0 &&\qquad\text{ in } A_{m}^0\cup E_{m}^0\text{ and}\\
				\tilde{f}_{m,\alpha} &= f\circ\pi_{m}^{m-1} &&\qquad\text{ in } D_{m}^{m-1}\cup F_{m}^{m-1}
		\end{alignat*}
	by construction, and therefore
		\begin{equation}
			\J_{A_{m}^0\cup E_{m}^0\cup D_{m}^{m-1}\cup F_{m}^{m-1}}(\tilde{f}_{m,\alpha})\leq 	\gamma\cdot a_m^3.
		\end{equation}
	
	Along the same line of reasoning, which is possible because of (\ref{eq:O-1-3-4}), we find 
		\begin{equation}\label{eq:O-1-4-16}
			\J_{E_{m}^j}(\tilde{f}_{m,\alpha}) \leq \gamma\cdot a_m^3,\qquad j=1,\dotsc,m-1
		\end{equation}
	and
		\begin{equation}\label{eq:O-1-4-17}
			\J_{F_{m}^j}(\tilde{f}_{m,\alpha}) \leq \gamma\cdot a_m^3, \qquad j=0,\dotsc,m-2.\medskip
		\end{equation}
	
	We now proceed with estimates on \(A_{m}^j\) \((j\neq0)\), and note that \(D_{m}^j\) \((j\neq m-1)\) can be treated analogously by symmetry. While using the Cube-Lemma (\Cref{lem:cube}) in the construction's background, we changed the original lift \(n\colon\Omega\to S^2\) of \(R\) into a Lipschitz mapping \(n_{m,\alpha}\colon \Omega\to S^2\) in order to get the new \(R_{m,\alpha}=F\circ n_{m,\alpha}\) having \(m\) singularities. The deformation part of the Cosserat energy on \(A_m^j\) (\(j\neq 0\)) thus is bounded once again by \(\gamma\cdot a_m^3\) with the same argument as for estimate (\ref{eq:O-1-4-14}). Additionally, for the micro-rotational part we have the following, notating \(\tilde{R}_{m,\alpha} = R_{m,\alpha}\circ\pi_{m}^j = F\circ n_{m,\alpha}\circ\pi_{m}^j=F\circ\tilde{n}_{m,\alpha}\) and with the fact that \(\tilde{n}_{m,\alpha}\) is constant in \((x_1,x_2,2ja_m-x_3)\)-direction. It holds
		\begin{align*}
			\frac{\partial\tilde{n}_{m,\alpha}}{\partial x_3} = \frac{x_1}{x_3-2ja_m}\frac{\partial\tilde{n}_{m,\alpha}}{\partial x_1} &+ \frac{x_2}{x_3-2ja_m}\frac{\partial\tilde{n}_{m,\alpha}}{\partial x_2}\\[5pt]
				\shortintertext{and thus}
			\abs{\D \tilde{n}_{m,\alpha}(x)}^2 = \left(1+\left(\tfrac{x_1}{x_3-2ja_m}\right)^2\right) \abs{\frac{\partial\tilde{n}_{m,\alpha}}{\partial x_1}}^2 &+ \left(1+\left(\tfrac{x_2}{x_3-2ja_m}\right)^2\right) \abs{\frac{\partial\tilde{n}_{m,\alpha}}{\partial x_2}}^2 \\[5pt]
				&+ 2\frac{x_1x_2}{(x_3-2ja_m)^2}\frac{\partial\tilde{n}_{m,\alpha}}{\partial x_1}\frac{\partial\tilde{n}_{m,\alpha}}{\partial x_2}.
		\end{align*}
	
	 As \(x_1^2+x_2^2\leq(x_3-2ja_m)^2\) in this regime, it directly follows that
		\begin{align*}
			\abs{\D(n_{m,\alpha}\circ\pi_{m}^j)(x)}^2 &\leq 3\cdot\left(\abs{\frac{\partial\tilde{n}_{m,\alpha}}{\partial x_1}}^2 + \abs{\frac{\partial\tilde{n}_{m,\alpha}}{\partial x_2}}^2\right)\\
				&= 3\cdot\left(\frac{a_m}{x_3-2ja_m}\right)^2\cdot\abs{\D n_{m,\alpha\vert B^2_{\frac{\alpha}{2}}\times\{(-1+2j)a_m\}}\left(\pi_{m}^j(x)\right)}^2\\
			&\leq 6\cdot\frac{8\alpha^4}{\left(\alpha^4+(\pi_{m}^j(x))_1^2+(\pi_{m}^j(x))_2^2\right)^2},
		\end{align*}
	using (\ref{eq:O-1-3-5}) and \(R_{m,\alpha}=F\circ n_{m,\alpha}\) in combination with the fact, that the covering map \(F\) is homothetic. With the transformation
		\begin{align*}
			\xi &=\sqrt{x_1^2+x_2^2+(x_3-2ja_m)^2},\\
				y_1 &= \left(\pi_{m}^j(x)\right)_1 = \frac{a_m}{\abs{x_3-2ja_m}}\cdot x_1 = \frac{a_m}{\sqrt{\xi^2-x_1^2-x_2^2}}\cdot x_1,\\
			y_2 &= \left(\pi_{m}^j(x)\right)_2 = \frac{a_m}{\abs{x_3-2ja_m}}\cdot x_2 = \frac{a_m}{\sqrt{\xi^2-x_1^2-x_2^2}}\cdot x_2,\\
			\eta^2 &= y_1^2+y_2^2
		\end{align*}
	and
		\begin{equation*}
			\dif x_1 \, \dif x_2 \, \dif x_3 = \frac{\xi^2 a_m\eta}{(a_m^2+\eta^2)^{\frac{3}{2}}} \, \dif\eta \, \dif\vartheta \, \dif\xi,
		\end{equation*}	
	we finally get 						
		\begin{align}\label{eq:O-1-4-18}
			&\J_{A_{m}^j}(\tilde{f}_{m,\alpha}) = \int\limits_{A_{m}^j} \abs{(R_{m,\alpha}\circ\pi_{m}^j)^T\cdot\D(\varphi\circ\pi_{m}^j)-I_3}^2 + \abs{\D(R_{m,\alpha}\circ\pi_{m}^j)}^2 \dif x \nonumber\\
				&\leq \gamma\cdot a_m^3 + \int\limits_{A_{m}^j} \abs{\D(F\circ(n_{m,\alpha}\circ\pi_{m}^j))}^2 \dif x = \gamma\cdot a_m^3 + 8\int\limits_{A_{m}^j} \abs{\D(n_{m\alpha}\circ\pi_{m}^j)}^2 \, \dif x \nonumber\\
			&\leq \gamma\cdot a_m^3 + 8\cdot 48\alpha^4\cdot \int\limits_0^{2\pi} \int\limits_0^{\frac{\alpha}{2}} \int\limits_{a_m}^{\sqrt{\eta^2+a_m^2}} \frac{1}{(\alpha^4+\eta^2)^2}\cdot\frac{\xi^2 a_m\eta}{(a_m^2+\eta^2)^{\frac{3}{2}}}\quad \dif\xi \, \dif\eta \, \dif\vartheta	\nonumber\\
				&=	\gamma\cdot a_m^3 + 8\cdot 32\pi\alpha^4\cdot \int\limits_0^{\frac{\alpha}{2}} \frac{a_m\eta}{(\alpha^4+\eta^2)^2}\cdot\frac{(a_m^2+\eta^2)^{\frac{3}{2}}-a_m^3}{(a_m^2+\eta^2)^{\frac{3}{2}}} \quad \dif\eta\nonumber\\
			&=	\gamma\cdot a_m^3 + 8\cdot 32\pi\alpha^4\cdot \int\limits_0^{\frac{\alpha}{2}} \frac{a_m\eta}{(\alpha^4+\eta^2)^2}\cdot\Bigg(1-\frac{1}{\left(1+(\frac{\eta}{a_m})^2\right)^{3/2}}\Bigg) \quad \dif\eta\nonumber\\
				&\leq \gamma\cdot a_m^3 + \gamma\cdot \alpha^4 \cdot \int\limits_0^{\frac{\alpha}{2}} \frac{a_m\eta}{(\alpha^4+\eta^2)^2}\cdot \frac{\eta^2}{a_m^2} \dif\eta,			
		\end{align}
	since \(\eta\leq\tfrac{\alpha}{2}< a_m\). Hence for \(j=1,\dotsc,m-1\) (\ref{eq:O-1-4-18}) becomes				
		\begin{align}\label{eq:O-1-4-19}
			\J_{A_{m}^j}(\tilde{f}_{m,\alpha}) &\leq \gamma\cdot a_m^3 + \gamma\alpha^4\cdot\frac{1}{a_m} \int\limits_0^{\frac{\alpha}{2}} \frac{\eta^3}{(\alpha^4+\eta^2)^2} \dif\eta \nonumber\\
				&\leq \gamma\cdot a_m^3 + \gamma\cdot\frac{\alpha^4}{a_m}\cdot\ln\left(1+\frac{1}{4\alpha^2}\right)\nonumber\\
			&\leq \gamma\cdot (a_m^3 + a_m).	
		\end{align}
	As mentioned above we also find 
		\begin{equation}\label{eq:O-1-4-20}
			\J_{D_{m}^j}(\tilde{f}_{m,\alpha})\leq \gamma\cdot(a_m^3+a_m)
		\end{equation}
	for \(j=0,\dotsc,m-2\) by symmetry.
	
	Combining (\ref{eq:O-1-4-13}) -- (\ref{eq:O-1-4-17}), (\ref{eq:O-1-4-19}) and (\ref{eq:O-1-4-20}), we have
		\begin{alignat*}{2}
			&\J_{K_m}(\tilde{f}_{m,\alpha}) = \sum_{j=0}^{m-1} \J_{\cmj}(\tilde{f}_{m,\alpha}) && \\
				&= \sum_{j=0}^{m-1} \left(\J_{B^3_{a_m}(c_j)}(\tilde{f}_{m,\alpha}) + \J_{G_{m}^j}(\tilde{f}_{m,\alpha})\right) && + \sum_{j=1}^{m-1} \left(\J_{A_{m}^j}(\tilde{f}_{m,\alpha}) + \J_{E_{m}^j}(\tilde{f}_{m,\alpha})\right)\\
			&+ \sum_{j=0}^{m-2} \left(\J_{D_{m}^j}(\tilde{f}_{m,\alpha}) + \J_{F_{m}^j}(\tilde{f}_{m,\alpha})\right) &&+ \J_{A_{m}^0\cup E_{m}^0\cup D_{m}^{m-1}\cup F_{m}^{m-1}}(\tilde{f}_{m,\alpha}) \\
				&\leq 2m\cdot 	a_m\cdot\left(1+\frac{\alpha^2}{a_m^2}\right)^2\left(64\pi+\varepsilon\right) &&+ \gamma\cdot(a_m^3+a_m)\\
			&\qquad\qquad\longrightarrow 2m\cdot a_m(64\pi+\varepsilon) &&+ \gamma\cdot(a_m^3+a_m), \text{ with }\alpha\searrow 0.
		\end{alignat*}
	In other words, for any \(\varepsilon>0 \text{ and any } m\geq m_0\), there is a mapping \(\tilde{f}_m\) with
		\begin{align*}
			\J_{K_m}(\tilde{f}_m) &\leq 2m\cdot a_m(64\pi+\varepsilon)+ \gamma\cdot(a_m^3+a_m)+ \varepsilon\\
				&\qquad\longrightarrow d\cdot 64\pi + (d+1)\varepsilon, \text{ with } m\to\infty.
		\end{align*}
	Hence, for each \(\varepsilon>0\) we have a number \(m\, (=m(\varepsilon))\geq m_0\) and a mapping  \(\tilde{f}_m\) that fulfils 
		\begin{equation*}
			\J_{K_m}(\tilde{f}_m)\leq 64\pi\cdot d + (d+2)\varepsilon.
		\end{equation*}
	Thus for any sequence \((\varepsilon_m)_{m\in\mathbb{N}}\) with \(\varepsilon_m \searrow 0\), we have constructed a sequence of Lipschitz mappings \(\tilde{f}_m\in H^1(\Omega,\mfk)\), such that
		\begin{enumerate}[label=(\roman*)]
			\item \begin{equation*}
				\tilde{f}_m =(\tilde{\varphi}_m,\tilde{R}_m)= \begin{cases}
					(\varphi,R)=f &\text{ in }\Omega\setminus K_{m},\\
				(\varphi\circ\pi_{m}^j, R_{m}\circ\pi_{m}^j) & \text{ in }K_{m},
				\end{cases}
				\end{equation*} where \(K_m\to[P,N], m\to\infty,\) in Hausdorff-distance,
			\item \(\sing(\tilde{f}_m) = \{P=c_0,c_1,\dotsc,c_{m-2},c_{m-1}=N\}\) with
				\begin{equation*}
					\begin{cases}
					\deg_{c_j}(\tilde{\varphi}_m)=0,	& j=0,\dotsc,m-1,\\
					\deg_{c_j}(\tilde{R}_m)=0,	& j=1,\dotsc,m-2,\\
					\deg_P(\tilde{R}_m)=1=\deg_N(\tilde{R}_m),
					\end{cases}
				\end{equation*}
			\item \begin{equation*}
				\limsup\limits_{m\to\infty} \J_{\Omega}(\tilde{f}_m)\leq\J_{\Omega}(f)+64\pi\cdot\abs{P-N}.
				\end{equation*}			
		\end{enumerate}
	
	\textbf{Step 3} (Approximation). Finally, we need suitable approximation arguments to achieve smoothness except in \(P,N\) for each \(\tilde{f}_m\) without affecting the Cosserat energy estimates. Also, \((P,N)\) is not yet a dipole for \(\tilde{R}_m\) according to \Cref{def:dipol}. Luckily, we are able to use several methods developed in \cite{bethuelapprox}. During the construction in Step 1, we changed the original smooth mapping \(f\) in the cuboid \(K_{m}\) only. Since \(H^1_{\tilde{\varphi}_m}(K_m,\mathbb{R}^3)\cap C^{\infty}(K_m,\mathbb{R}^3)\) is dense in \(H^1_{\tilde{\varphi}_m}(K_m,\mathbb{R}^3)\), it is possible to approximate the changed deformation component \(\tilde{\varphi}_{m\vert K_{m}}\) in \(H^1\)-topology with mappings \[\varphi_{m,s}\in C^{\infty}(K_{m},\mathbb{R}^3),\quad s\in\mathbb{N},\] while the boundary values remain 
	\(\varphi_{m,s} = \tilde{\varphi}_m=\varphi \text{ in }\partial  K_{m}\).
	Replacing \(\varphi_{m,s}\) by a subsequence, we may assume
	\(\varphi_{m,s} \to \tilde{\varphi}_m\), \(s\to\infty\), pointwise almost everywhere.\\
	
	For the micro-rotational component \(\tilde{R}_{m\vert K_{m}} = (R_{m}\circ\pi_{m}^j)\in H^1(K_{m},\mathcal{S})\) we first note again that \(\tilde{R}_{m\vert K_{m}}= F\circ n_{m} \circ\pi_{m}^j = F\circ\tilde{n}_m,\) with \(\tilde{n}_m\in H^1(K_{m}, S^2)\), \[\deg_P(\tilde{n}_m) =2k+1=-\deg_N(\tilde{n}_m),\text{ for some } k\in\mathbb{Z}, \text{ and } \deg_{c_j}(\tilde{n}_m)=0,\] for each \(j=1,\dotsc,m-2,\) having smooth boundary values \(\tilde{n}_{m\vert \partial K_{m}}=n_{\vert \partial K_{m}}\in C^{\infty}(\partial K_{m},S^2)\) due to the underlying construction making use of the smooth lift \(n\) of the original smooth \(R\).			
	Applying \cite[][Theorem 2 bis]{bethuelapprox}, since \(S^2\) is a compact manifold without boundary and \(K_{m}\) is dividable into cubes (cf. "cubeulation" in \cite{bethuelapprox}), there exists a sequence \[n_{m,t}\in H^1(K_{m},S^2)\cap C^{\infty}(K_{m}\setminus\{P,N,c_1,\dotsc,c_{m-2}\}, S^2)\text{ with }\]
		\begin{enumerate}[label=(\roman*)]
			\item \(n_{m,t}\to\tilde{n}_m,\) \(t\to\infty,\quad\)  in \(H^1(K_{m},S^2)\),
				\item \(\deg_P(n_{m,t})= 2k+1 = -\deg_N(n_{m,t})\) and\\ \(\deg_{c_j}(n_{m,t})=0\), \(j=1,\dotsc,m-2\),
			\item \(n_{m,t\vert \partial K_{m}}=\tilde{n}_{m\vert \partial K_{m}} = n_{\vert \partial K_{m}}.\)
		\end{enumerate}	
	(Readers interested in details of \cite[][Theorem 2 bis]{bethuelapprox} should pay attention to a typing error there, which can easily be spotted when comparing with \cite[Theorem 1 bis,][]{bethuelapprox}. The correct (and fulfilled) assumption for our case is that \(\tilde{n}_m\) restricted to \(\partial K_m\) is smooth in \(\partial K_m\), not in the whole of \(K_m\).)
	
	Now we can use the technique from the proof of \cite[][Lemma 1 bis]{bethuelapprox} to get rid of those singularities \(c_1,\dotsc,c_{m-2}\), in which the homotopy class of \(n_{m,t}\) is trivial. For \(n_{m,t}\) in \(Q_{m}\coloneqq\bigcup\limits_{j=1}^{m-2} \cmj\), each \(n_{m,t\vert Q_{m}}\) (subject to their own boundary values \(g\coloneqq n_{m,t\vert \partial Q_{m}}\)) can be approximated in \(H^1\)-topology by mappings \[n_{m,t,s}\in H^1_g(Q_{m},S^2)\cap C^{\infty}(Q_{m}, S^2),\] that agree with \(n_{m,t}\) outside of \(\bigcup\limits_{j=1}^{m-2} B^3_{1/s}(c_j)\).\medskip
	
	\noindent That is why for each \(\tilde{n}_m\colon K_{m}\to S^2\), there exists a sequence of mappings \[n_{m,t,s}\in H^1(K_{m},S^2)\cap C^{\infty}(K_{m}\setminus\{P,N\}, S^2),\quad s\in\mathbb{N}, \] by defining \[ n_{m,t,s} = \begin{cases}
		n_{m,t,s} & \text{ in } Q_{m},\\
		n_{m,t} & \text{ in } C_{m}^0\cup C_{m}^{m-1},
	\end{cases} \] with \(n_{m,t,s}\to\tilde{n}_m\) in \(H^1(K_{m}, S^2)\), \(\deg_P(n_{m,t,s})=2k+1=-\deg_N(n_{m,t,s})\) and \(n_{m,t,s\vert \partial K_{m}} = \tilde{n}_{m\vert \partial K_{m}} = n_{\vert \partial K_{m}}\). \\
	
	Finally we project everything back from \(S^2\) to \(\mathcal{S}\). For each of the mappings \(\tilde{R}_{m\vert K_m}=F\circ\tilde{n}_m\colon K_{m}\to\mathcal{S}\), we thus have a sequence of mappings \[R_{m,t,s}\coloneqq F\circ n_{m,t,s}\in H^1(K_{m},\mathcal{S})\cap C^{\infty}(K_{m}\setminus\{P,N\},\mathcal{S}),\] which approximate \(\tilde{R}_{m\vert K_m}\) in \(H^1\)-topology, because (after passing to an a.e.-pointwise convergent subsequence) it holds that	
		\begin{align*}
			&\int\limits_{K_{m}} \abs{R_{m,t,s}-\tilde{R}_m}^2 \dif x =\int\limits_{K_{m}} \abs{F\circ n_{m,t,s}-F\circ\tilde{n}_m}^2 \dif x\\
				&= \int\limits_{K_{m}} \abs{2\left[n_{m,t,s}\otimes n_{m,t,s} -\tilde{n}_m\otimes \tilde{n}_m \right]}^2\dif x\\
			&\leq  \sum\limits_{i,j=1}^3 \quad \int\limits_{K_{m}} \big[ \underbrace{\abs{n_{m,t,s}^{i}-\tilde{n}_m^{i}}}_{\to 0} \cdot \underbrace{\abs{n_{m,t,s}^{j}+\tilde{n}_m^{j}}}_{\leq 2} + \underbrace{\abs{n_{m,t,s}^{i}+\tilde{n}_m^{i}}}_{\leq 2} \cdot \underbrace{\abs{n_{m,t,s}^{j}-\tilde{n}_m^{j}}}_{\to 0}\big]^2 \dif x\\
				&\longrightarrow 0,\, s\to\infty\\[5pt]
		\shortintertext{and}
			&\int\limits_{K_{m}} \abs{\D R_{m,t,s}-\D\tilde{R}_m}^2 \dif x =\int\limits_{K_{m}} \abs{\D F\circ n_{m,t,s}\cdot \D n_{m,t,s} - \D F\circ\tilde{n}_m\cdot \D\tilde{n}_m}^2 \dif x\\
				&\leq 2\int\limits_{K_{m}} \abs{\D F\circ n_{m,t,s}\cdot \D n_{m,t,s} - \D F\circ\tilde{n}_m\cdot \D n_{m,t,s}}^2\\
			&\qquad\qquad+ \abs{\D F\circ \tilde{n}_m\cdot \D n_{m,t,s} - \D F\circ\tilde{n}_m\cdot \D\tilde{n}_m}^2 \dif x\\[5pt]
				&\leq 2\int\limits_{K_{m}} \abs{ \D F\circ n_{m,t,s}- \D F\circ\tilde{n}_m }^2\cdot\abs{ \D n_{m,t,s} }^2_{\infty} + \abs{ \D F\circ \tilde{n}_m }^2_{\infty}\cdot \abs{\D n_{m,t,s} - \D\tilde{n}_m}^2 \dif x\\
			&\longrightarrow 0,\, s\to\infty.
		\end{align*}
	
	Summarizing, for any sequence \((\varepsilon_m)_{m\in\mathbb{N}}\) with \(\varepsilon_m\searrow 0\), there exists a sequence of Sobolev mappings \((\tilde{f}_m)_m \text{ with } \) \[\J_{\Omega}(\tilde{f}_m)\leq \J_{\Omega}(f) + 64\pi\cdot d + (d+2)\varepsilon_m,\] as well as another sequence of mappings \((f_{m,t,s})_s\), \[f_{m,t,s}\in H^1(\Omega,\mfk)\cap C^{\infty}(\Omega\setminus\{P,N\},\mfk),\]  \[f_{m,t,s}\coloneqq \begin{cases}
		(\varphi,R)=f, & \text{ in }\Omega\setminus K_{m}\\
		(\varphi_{m,t,s}, R_{m,t,s}), & \text{ in } K_{m} 
	\end{cases}\] with \(f_{m,t,s}\to\tilde{f}_m\), \(s\to\infty\), in \(H^1(\Omega,\mfk)\). Moreover, \(R_{m,t,s}\) has a dipole \((P,N)\).
	
	Hence for each \(\varepsilon_m\), by dominated convergence, we get the existence of a mapping \(f_m\) with the desired properties and
		\begin{align*}
			\J_{K_{m}}(f_m) &\leq \J_{K_{m}}(\tilde{f}_m) + \varepsilon_m\\
				&\leq 64\pi\cdot d + (d+3)\varepsilon_m,
		\end{align*}
	meaning we have found the sequence \(\left(f_m\right)_{m\in\mathbb{N}}\) of mappings, which are smooth except for a dipole in the micro-rotation and whose Cosserat energy fulfils \[\lim\limits_{m\to\infty} \J_{\Omega}(f_m)\leq \J_{\Omega}(f) + 64\pi\cdot\abs{P-N}.\]
\end{proof}

Having finished the proof of \Cref{thm:dipoles}, it remains to prove the Cube Lemma for the Cosserat energy.

\begin{proof}[Proof of \Cref{lem:cube}.]
	Since the set \(\partial C_{\nu}\) is bounded, simply connected and locally path-connected, the given Lipschitz mapping \(R\colon\partial C_{\nu}\to\mathcal{S}\) can be lifted, which means for the covering map \(F\) of \(\mathcal{S}\) (\(F\colon S^2\to\mathcal{S}, q\mapsto 2q\otimes q-I_3\)), there exist exactly two Lipschitz mappings
		\begin{equation*}
			n_i\colon\partial C_{\nu}\to S^2 \qquad\text{ with }\qquad R=F\circ n_i,
		\end{equation*} \(i=1,2\) and \(n_1=-n_2\) (cf. \cref{sec1:introduction}). We choose one of these mappings and keep it fixed (\(n:=n_1\)). Additionally, there is a number \(\deg(n)=c_0\in \mathbb{Z}\), such that \(d_0 =\deg(R)= c_0\mod 2\). Moreover, \(F\) is homothetic, i.e. for any tangent vector \(V\in T_p(\partial C_{\nu})\) we have
		\begin{equation}\label{eq:O-1-3-7}
			\abs{\D R_p(V)}^2 = \abs{\D F_{n(p)}(\D n_p(V))}^2 = 8\cdot\abs{\D n_p(V)}^2.
		\end{equation}
	
	We perform the modifications used in \cite{brezis1983large} on our \(n\), to the following effect. For each \(\varepsilon>0\) there is a constant \(0<\alpha_0<\nu\), such that for each \(0<\alpha<\alpha_0\), there exists a Lipschitz mapping \(\tilde{n}\colon\partial C_{\nu}\to S^2\), which fulfils 
		\begin{align}
			\deg(\tilde{n}) &= c_0+1,\nonumber\\
				\tilde{n} &= n \text{ on } \partial C_{\nu}\setminus\left(\ball\right)\nonumber\\
		\shortintertext{and}
			\int\limits_{\ball}  \abs{\D\tilde{n}}^2 \,\dif\mathcal{H}^2 &= 8\pi + O(\alpha^2) \qquad\text{ with }\alpha\searrow 0,\label{eq:O-1-3-10}\\
		\shortintertext{as well as}
			\abs{\D\tilde{n}} &\leq const \qquad \text{ in }\ring \label{eq:O-1-3-11}\\
		\shortintertext{and}
			\abs{\D\tilde{n}(x,y,0)}^2 &= \frac{8\alpha^4}{(\alpha^4+x^2+y^2)^2} \qquad \text{ in }\sball \label{eq:O-1-3-12}.
		\end{align}
	That is to say,  we insert a singularity of degree \(+1\) into \(n\), without changing \(n\) outside of the disc \(\ball\). For doing so, we only need a controlled amount of (Dirichlet) energy. Now we define \(\tilde{R}\coloneqq F\circ\tilde{n}\), such that
		\begin{equation*}
			\tilde{R} = F\circ\tilde{n} = F\circ n = R \text{ on } \partial C_{\nu}\setminus\left(\ball\right)
		\end{equation*}
	and
		\begin{equation*}
			\deg(\tilde{R}) = \deg(\tilde{n}) \mod 2 = c_0+1 \mod 2 = d_0+1 \mod 2.
		\end{equation*}
	The properties (\ref{eq:O-1-3-4}) and (\ref{eq:O-1-3-5}) follow directly from (\ref{eq:O-1-3-11}) and (\ref{eq:O-1-3-12}) respectively, because of (\ref{eq:O-1-3-7}).	For the other part of the Cosserat energy of \(\tilde{f}\), coming from the deformation, we easily see
	\(\abs{\tilde{R}^T\cdot \D\varphi - I_3}^2 = O(1), \alpha\searrow 0\), as \(\tilde{R}\in \mathcal{S}\subset SO(3)\) and  \(\varphi\) is Lipschitz by assumption.
	
	Using (\ref{eq:O-1-3-7}) and (\ref{eq:O-1-3-10}), we therefore find
		\begin{align*}
			\int\limits_{\ball} 2\abs{\tilde{R}^{T}\cdot \D\varphi - I_3}^2 + \abs{\D\tilde{R}}^2 \;\dif\mathcal{H}^2
				&= \int\limits_{\ball} 2\abs{\tilde{R}^{T}\cdot \D\varphi - I_3}^2 + 8\cdot\abs{\D\tilde{n}}^2 \;\dif\mathcal{H}^2\\
			&= 64\pi + O(\alpha^2) \qquad\text{ with } \alpha\searrow 0.
		\end{align*}		
	Thus we can choose \(\alpha_0\ll 1\) sufficiently small, such that \[\int_{\ball} 2\abs{\tilde{R}^{T}\cdot \D\varphi - I_3}^2 + \abs{\D\tilde{R}}^2 \;\dif\mathcal{H}^2 < 64\pi + \varepsilon \] holds for every \(\alpha<\alpha_0\).
	
\end{proof}

\section{Prescribing the number of singular points in a Cosserat-elastic solid (proof)}\label{sec:resmin}
Now that we have provided a tool for constructing dipoles with a controlled amount of Cosserat energy, we can use it to prove \Cref{thm:HL}.
\begin{proof}[Proof of \Cref{thm:HL}]
	We use a combination of ideas from \cite{brezis1989s} and \cite{hardt1986remark}, adapted to the restricted Cosserat problem (\ref{eq: restricted problem}). First, we define the desired smooth boundary conditions following an idea from \cite[][Section II.4,]{brezis1989s}. Then, we show that each restricted minimizer of the Cosserat energy in the corresponding class must at least have a given number of point singularities, just as it was done in \cite{hardt1986remark} for harmonic mappings.\medskip
	
	For \(N\in\mathbb{N}\), we define numbers \(\lambda_i \coloneqq\frac{i}{2N}\), \(i=1,\dotsc,N\), points \[\xi_i = (0, \sqrt{1-\lambda_i^2}, \lambda_i) \qquad\text{ and }\qquad
	\eta_i = (0 , -\sqrt{1-\lambda_i^2},  -\lambda_i),\] as well as pairs of points
		\[\begin{matrix}
			P_{i,\varepsilon}^{+}=(1-\varepsilon)\xi_i \\ N_{i,\varepsilon}^{+}=(1+\varepsilon)\xi_i
		\end{matrix} \qquad\text{ and }\qquad \begin{matrix}
			P_{i,\varepsilon}^{-}=(1+\varepsilon)\eta_i \\ N_{i,\varepsilon}^{-}=(1-\varepsilon)\eta_i,
		\end{matrix}\]
	for each \(\varepsilon>0\). By \(Z_{i,\varepsilon}^{+}\) and \(Z_{i,\varepsilon}^{-}\), we denote the \(\varepsilon\)-tubular neighbourhood of the line segment \(\left[P_{i,\varepsilon}^{+} , N_{i,\varepsilon}^{+}\right]\) and \(\left[P_{i,\varepsilon}^{-} , N_{i,\varepsilon}^{-}\right]\) respectively.
	
	Without loss of generality, let \(\varepsilon\ll 1\) such that
		\begin{equation}\label{eq:O-1-5-01}
			\abs{z-\hat{z}}\geq \frac{1}{4N} \quad \text{for all}\, (x,y,z)\in Z_{i,\varepsilon}^{+},\, (\hat{x},\hat{y},\hat{z})\in Z_{j,\varepsilon}^{+}\, \text{with } i\neq j,
		\end{equation}
	meaning the neighbourhoods \(Z_{i,\varepsilon}^{+}\) are separated by circular slices with a height of at least \(\tfrac{1}{4N}\). Moreover, let \(g=(\vartheta,S)\in C^{\infty}\left(B^3_2,\mfk\right)\) be given as
		\[\vartheta(x,y,z) \equiv \begin{pmatrix}
			-x\\-y\\z \end{pmatrix} \quad \text{ and } \quad S(x,y,z) \equiv \begin{pmatrix}
			-1 & 0  & 0 \\
			0 & -1 & 0 \\
			0 & 0  & 1 \end{pmatrix}.\]
	At first, \Cref{thm:dipoles} gives us a mapping \(\tilde{g} = (\tilde{\vartheta}, \tilde{S})\) in \[H^1(B^3_2,\mfk)\cap C^{\infty}\left(B^3_2\setminus\left(\bigcup\limits_{i=1}^N\left\{P_{i,\varepsilon}^{+},N_{i,\varepsilon}^{+},P_{i,\varepsilon}^{-},N_{i,\varepsilon}^{-} \right\}\right),\mfk\right),\] which agrees with \(g\) outside of \(Z\coloneqq\bigcup\limits_{i=1}^N \left(Z_{i,\varepsilon}^{+}\cup Z_{i,\varepsilon}^{-}\right)\), with \(\tilde{S}\) having only the dipoles \((P_{i,\varepsilon}^{\pm},N_{i,\varepsilon}^{\pm})\) as singularities,
		\begin{align}
			\deg_{P_{i,\varepsilon}^{+}}(\tilde{S}) &= 1 = \deg_{N_{i,\varepsilon}^{+}}(\tilde{S}),\nonumber\\
				\deg_{P_{i,\varepsilon}^{-}}(\tilde{S}) &= 1 = \deg_{N_{i,\varepsilon}^{-}}(\tilde{S})\nonumber\\
		\shortintertext{and}
			\J_{B^3}(\tilde{g}) < 64\pi\cdot &2\varepsilon \cdot 2N + \varepsilon = (64\pi\cdot 4N + 1)\cdot \varepsilon.\label{eq:O-1-5-04}
		\end{align}
	Now, the boundary values we prescribe are \[g_0 = (\varphi_0,R_0) \coloneqq \tilde{g}_{\vert \partial B^3} \in C^{\infty}(\partial B^3, \mfk),\] with \(\deg(R_0)=0\) and  \[\tilde{g}_{\vert \overline{B^3}\setminus Z} \equiv (\vartheta,S).\]
	
	Due to (\ref{eq:O-1-5-04}), combined with (\ref{eq:O-1-5-01}), it is possible to choose  $\varepsilon\ll 1 $ such that
		\begin{equation}
			\J_{B^3}(\tilde{g}) < \frac{\pi}{N}. \label{eq:O-1-5-06}
		\end{equation}
	
	\noindent For the rest of this proof \(\tilde{g}\) always means \(\tilde{g}_{\vert \overline{B^3}}\) and therefore, \(\tilde{g}\in H^1_{g_0}(B^3,\mfk)\). Thus, let \(f=(\varphi,R)\) be a restricted minimizer, i.e. a minimizer of the problem (\ref{eq: restricted problem}) in the (non-empty) class \(H^1_{g_0}(B^3,\mfk)\). Because of (\ref{eq:O-1-5-06}), we also have
		\begin{equation}
			\J_{B^3}(f) < \frac{\pi}{N}. \label{eq:O-1-5-07}
		\end{equation}
	
	From here, following \cite{hardt1986remark}, with (\ref{eq:O-1-5-01}) and (\ref{eq:O-1-5-07}), we find \(N+1\) numbers \(\mu_0,\dotsc,\mu_N\in\mathbb{R}\), \[0<\mu_0<\lambda_1<\mu_1<\lambda_2<\dots<\lambda_N<\mu_N<1,\]
	such that for each \(i=0,\dotsc,N\) the disc \(D_i\coloneqq\left\{(x,y,z)\in B^3: z=\mu_i\right\}\) satisfies \(\overline{D}_i\cap Z = \emptyset\) and
		\begin{equation}
			\int\limits_{\overline{D}_i} \abs{\D R}^2 \dif\mathcal{H}^2 < 4\pi. \label{eq:O-1-5-09}
		\end{equation} 
	Suppose (\ref{eq:O-1-5-09}) was not possible. Then we would get \[\J_{B^3}(f)\geq\frac{1}{4N}\cdot4\pi = \frac{\pi}{N},\] because of (\ref{eq:O-1-5-01}), contradicting (\ref{eq:O-1-5-07}). \smallskip		
	
	As we mentioned in the introduction, following the line of reasoning from \cite{gastel2019regularity} and carrying out similar arguments at the boundary, we get discreetness of the singular set for restricted Cosserat energy minimizers not only in the interior, but up to the boundary. In particular, singularities cannot accumulate at the boundary. This implies that also \(\sing(R)\subseteq\sing(f)\subset\overline{B^3}\) is discrete in \(\overline{B^3}\). Hence we can assume without loss of generality, that all isolated singularities of \(R\) are in \(\overline{B^3}\setminus\left(\bigcup\limits_{i=1}^N \overline{D_i}\right)\). Then for each compact subset \(C\) of \(\overline{B^3}\setminus\sing(R)\) and each \(i=0,\dots,N\), we have
		\begin{align*}
			\int\limits_{\mathcal{S}} \mathcal{H}^0 \left(\overline{D}_i\cap C\cap R^{-1}(p)\right)\,\dif\mu(p) &= \int\limits_{\overline{D}_i\cap C} \operatorname{Jac} (R_{\vert \overline{D}_i\cap C}) \,\dif \mathcal{H}^2\\
				\leq \int\limits_{\overline{D}_i} \operatorname{Jac} (R_{\vert \overline{D}_i\cap C}) \,\dif \mathcal{H}^2 &\leq \frac{1}{2} \int\limits_{\overline{D}_i} \abs{\D R}^2 \,\dif\mathcal{H}^2,
		\end{align*}
	where \(\mu\) denotes the natural Riemannian measure on \(\mathcal{S}\). From monotone convergence (for \(C\nearrow B^3\)), together with (\ref{eq:O-1-5-09}), we infer for each \(i=0,\dotsc,N\) 	
		\begin{align*}
			\mu(R(\overline{D}_i)) &\leq \int\limits_{\overline{D}_i} \mathcal{H}^0 \left(\overline{D}_i\cap R^{-1}(p)\right) \,\dif\mu(p) = \int\limits_{\overline{D}_i}  \operatorname{Jac} (R_{\vert \overline{D}_i}) \,\dif \mathcal{H}^2\\
				&\leq \frac{1}{2} \int\limits_{\overline{D}_i} \abs{\D R}^2 \,\dif\mathcal{H}^2 <  2\pi = \mu(\mathcal{S}).
		\end{align*}
	
	So each image \(R(\overline{D}_i)\) is a proper subset of \(\mathcal{S}\). On one hand, as \(R\) is (Hölder-) continuous on \(\overline{D}_i\) and \(\overline{D}_i\cap Z = \emptyset\), combined with  \(R_{\vert \partial B^3}=R_0=\tilde{S}_{\vert \partial B^3}\) and \[\tilde{S}_{\vert \overline{B^3}\setminus Z} = S = \begin{pmatrix}
		-1 & 0 & 0\\
		0 & -1 & 0\\
		0 & 0 & 1
	\end{pmatrix},\] it holds that \(R_{\vert \overline{D}_i}\) is homotopic to \(S\) relative to \(\partial D_i\), i.e.
		\begin{equation}
			R_{\vert \overline{D}_i} \simeq \begin{pmatrix}
			-1 & 0 & 0\\
			0 & -1 & 0\\
			0 & 0 & 1
			\end{pmatrix} \text{ rel. } \partial D_i. \label{eq:O-1-5-10} 		
		\end{equation}
	And on the other hand, we infer
		\begin{equation}
			\tilde{S}_{\vert \overline{D}_i} \equiv \begin{pmatrix}
			-1 & 0 & 0\\
			0 & -1 & 0\\
			0 & 0 & 1
			\end{pmatrix}. \label{eq:O-1-5-11}
		\end{equation}
	
	Finally, for each \(i=1,\dotsc,N\), we consider the slices \[\Omega_i\coloneqq\left\{(x,y,z)\in B^3: \mu_{i-1}<z<\mu_i\right\}.\] Then \(\partial\Omega_i=D_{i-1}\cup D_i\cup \left(\partial B^3\cap\overline{\Omega}_i\right)\) is homeomorphic to \(S^2\). Since the only singularity of \(\tilde{S}\) in \(\Omega_i\) is the point \(P_{i,\varepsilon}^{+}\) by construction, we get \[\deg(\tilde{S}_{\vert \partial\Omega_i}) = 1.\] Moreover,  (\ref{eq:O-1-5-10}), (\ref{eq:O-1-5-11}) and \(R_{\vert \partial B^3} = \tilde{S}_{\vert \partial B^3}\) imply  \(R_{\vert \partial\Omega_i} \simeq \tilde{S}_{\vert \partial\Omega_i}\). Hence, due to homotopy invariance of \(\deg(\cdot)\), we also find
	\[\deg\left(R_{\vert \partial\Omega_i}\right) = \deg(\tilde{S}_{\vert \partial\Omega_i}) =  1.\] This proves that \(R\) must have at least one singularity in \(\Omega_i\) for each \(i=1,\dotsc,N\).
	
\end{proof}

\bmhead{Acknowledgments}

The author's work on this subject is part of a project funded by the Deutsche Forschungsgemeinschaft (DFG, German Research Foundation), project ID 441380936. We appreciate all the encouraging discussions with members of priority programm (SPP) 2256, to which this project belongs.

\bmhead{Data availability statement}

Data sharing is not applicable to this article as no datasets were generated or analysed during the current study.

\bibliography{huesken-bibliography} 

\end{document}